\overfullrule=0pt
\centerline {\bf Existence, uniqueness, localization and minimization property of positive solutions for non-local problems} 
\centerline {\bf involving discontinuous Kirchhoff functions}\par
\bigskip
\bigskip
\centerline {BIAGIO RICCERI}\par
\bigskip
\bigskip
{\bf Abstract:} Let $\Omega\subset {\bf R}^n$ be a smooth bounded domain. In this paper, we prove a result of which the following
is a by-product: Let $q\in ]0,1[$, $\alpha\in L^{\infty}(\Omega)$, with $\alpha>0$, and $k\in {\bf N}$. Then, the problem
$$\cases {-\tan\left(\int_{\Omega}|\nabla u(x)|^2dx\right)\Delta u=
\alpha(x)u^q
 & in
$\Omega$\cr & \cr u>0 & in $\Omega$\cr & \cr
u=0 & on
$\partial \Omega$
\cr & \cr
(k-1)\pi<\int_{\Omega}|\nabla u(x)|^2dx<(k-1)\pi+{{\pi}\over {2}}
\cr}$$
has a unique weak solution $\tilde u$ which is the unique global minimum in $H^1_0(\Omega)$ of the functional
$$u\to {{1}\over {2}}\tan\left (\int_{\Omega}|\nabla\tilde u(x)|^2dx\right)\int_{\Omega}|\nabla u(x)|^2dx-{{1}\over {q+1}}\int_{\Omega}\alpha(x)|u^+(x)|^{q+1}dx\ ,$$
where $u^+=\max\{0,u\}$.
\bigskip
\bigskip
\bigskip
\bigskip
{\bf 1. Introduction}\par
\bigskip
First, we stress that we have chosen the above long title just to summarize the main features and novelties of our results.\par
\smallskip
Throughout the sequel, $\Omega\subset {\bf R}^n$ is a smooth bounded domain and $K$ is a real-valued function defined in $[0,+\infty[$.
\par
\smallskip
 Given a function $\varphi:\Omega\times {\bf R}\to {\bf R}$, we are interested in the problem
$$\cases {-K\left(\int_{\Omega}|\nabla u(x)|^2dx\right)\Delta u=
\varphi(x,u)
 & in
$\Omega$\cr & \cr u>0 & in $\Omega$\cr & \cr
u=0 & on
$\partial \Omega$\ .\cr}$$
We define weak solution of this problem any $u\in H^1_0(\Omega)$, with $u>0$ in $\Omega$, such that, for every $v\in H^1_0(\Omega)$, the function $\varphi(\cdot,u(\cdot))v(\cdot)$ lies in $L^1(\Omega)$ and one has
 $$K\left(\int_{\Omega}|\nabla u(x)|^2dx\right)\int_{\Omega}\nabla u(x)\nabla v(x)dx
-\int_{\Omega}
\varphi(x,u(x))v(x)dx=0\ .$$
This is a Kirchhoff-type problem, $K$ being the Kirchhoff function. Unquestionably, it is among the most studied nonlinear problems of the last two dacades. For a lucid introduction to the subject jointly with the relevant bibliography, we refer to the recent [5].\par
\smallskip
Here is a most remarkable corollary of our main result.\par
\medskip
THEOREM 1.1. - {\it Assume that there exists an open interval $I\subseteq ]0,+\infty[$ such that the restriction of $K$ to $I$ is increasing
and $K(I)=]0,+\infty[$. \par
Then, for each $q\in ]0,1[$ and for each $\alpha\in L^{\infty}(\Omega)$, with $\alpha>0$, the problem
$$\cases {-K\left(\int_{\Omega}|\nabla u(x)|^2dx\right)\Delta u=\alpha(x) u^q
 & in
$\Omega$\cr & \cr u>0 & in $\Omega$\cr & \cr
u=0 & on
$\partial \Omega$\cr & \cr
\int_{\Omega}|\nabla u(x)|^2dx\in I\cr}$$
has a unique weak solution $\tilde u$ which is the unique global minimum in $H^1_0(\Omega)$ of the functional
$$u\to {{1}\over {2}}K\left(\int_{\Omega}|\nabla \tilde u(x)|^2dx\right)\int_{\Omega}|\nabla u(x)|^2dx-{{1}\over {q+1}}\int_{\Omega}\alpha(x)|u^+(x)|^{q+1}dx \ .\eqno {(1.1)}$$
Moreover, $\tilde u$ satisfies the inequality
$$\left (K\left(\int_{\Omega}|\nabla \tilde u(x)|^2dx\right)\right)^{{2}\over {1-q}}
\int_{\Omega}|\nabla \tilde u(x)|^2dx\leq \left [\left ( {{2}\over {q+1}}\right )^2\left ({{\hbox {\rm ess sup}_{\Omega}\alpha}\over 
{\lambda_1}}\right )^{q+1}
\right ]^{{1}\over {1-q}}\int_{\Omega}\alpha(x)dx\ ,$$
where
$$\lambda_1=\inf_{u\in H^1_0(\Omega)\setminus \{0\}}{{\int_{\Omega}|\nabla u(x)|^2dx}\over {\int_{\Omega}|u(x)|^2dx}}\ .$$}\par
\medskip
The main novelties of Theorem 1.1 are the lack of continuity of $K$ in $[0,+\infty[$, the localization of the solution given by 
 $\int_{\Omega}|\nabla u(x)|^2dx\in I$ and the property that $\tilde u$ minimizes the functional $(1.1)$ (which depends on $\tilde u$ itself).
\smallskip
Actually, as far as we know, the continuity of $K$
in $[0,+\infty[$ is an assumption present in each paper devoted to this subject, but [9]. More precisely, in [9], we assumed that, for some $r>0$, $K$ is continuous and increasing in $[0,r[$, with $\lim_{t\to r^-}\int_0^tK(s)ds=+\infty$.\par
\smallskip
So, for instance, given $q\in ]0,1[$, $\alpha\in L^{\infty}(\Omega)$, with $\alpha>0$, $k\in {\bf N}$,  $c>0$
and $s\in ]0,1[$, the problems
$$\cases {-\log\left(\int_{\Omega}|\nabla u(x)|^2dx\right)\Delta u=
\alpha(x)u^q
 & in
$\Omega$\cr & \cr u>0 & in $\Omega$\cr & \cr
u=0 & on
$\partial \Omega$
\cr & \cr
\int_{\Omega}|\nabla u(x)|^2dx>1
\cr}$$
$$\cases {-\tan\left(\int_{\Omega}|\nabla u(x)|^2dx\right)\Delta u=
\alpha(x)u^q
 & in
$\Omega$\cr & \cr u>0 & in $\Omega$\cr & \cr
u=0 & on
$\partial \Omega$
\cr & \cr
(k-1)\pi<\int_{\Omega}|\nabla u(x)|^2dx<(k-1)\pi+{{\pi}\over {2}}
\cr}$$
and
$$\cases {-\left|c-\int_{\Omega}|\nabla u(x)|^2dx\right|^{-s}\Delta u=\alpha(x)u^q
 & in
$\Omega$\cr & \cr u>0 & in $\Omega$\cr & \cr
u=0 & on
$\partial \Omega$
\cr & \cr
\int_{\Omega}|\nabla u(x)|^2dx<c\cr}$$
cannot be covered by any of the results known up to now. Thanks to Theorem 1.1, each of these problems has
a unique weak solution $\tilde u$ which minimizes the functional $(1.1)$.
\smallskip
But even when $I=]0,+\infty[$, Theorem 1.1 turns out to be new. Actually, the known results for the problem
$$\cases {-K\left(\int_{\Omega}|\nabla u(x)|^2dx\right)\Delta u=u^q
 & in
$\Omega$\cr & \cr u>0 & in $\Omega$\cr & \cr
u=0 & on
$\partial \Omega$\cr}$$
where $q\in ]0,1[$, require that $K$, besides being continuous in $[0,+\infty[$, is non-increasing ([1], Theorem 4; [3], Example 1) or 
$K(t)=at+b$, with $a, b>0$  ([10], Theorem 1.2; [11], Corollary 1.1), or with $a>0$ and $b\geq 0$ ([4], Theorem 5.4).
\bigskip
{\bf 2. Results}\par
\bigskip
If $X$ is a topological space, a function $f:X\to {\bf R}$ is said to be inf-compact (resp. sup-compact) provided that, for each $r\in {\bf R}$,
the set $f^{-1}(]-\infty,r])$ (resp. $f^{-1}([r,+\infty[)$) is compact.\par
\smallskip
We will obtain our main result via the following abstract theorem.\par
\medskip
THEOREM 2.1. - {\it Let $X$ be a topological space and let $\Phi:X\to {\bf R}$, with $\Phi^{-1}(0)\neq \emptyset$, and $J:X\to {\bf R}$ be two functions such that, for each $\lambda>0$, the function $\lambda \Phi-J$ is lower semicontinuous, inf-compact and
admits a unique global minimum. Moreover, assume that $J$ has no global maxima in $X$. Furthermore, let $I\subseteq ]0,+\infty[$
be an open interval and $\Psi:I\to {\bf R}$ be an increasing function such that $\Psi(I)=]0,+\infty[$.\par
Then, there exists a unique $\tilde x\in X$  such that $\Phi(\tilde x)\in I$ and
$$\Psi(\Phi(\tilde x))\Phi(\tilde x)-J(\tilde x)=\inf_{x\in X}(\Psi(\Phi(\tilde x))\Phi(x)-J(x))\ .$$}\par
\smallskip
PROOF. Clearly, the function $\Psi^{-1}$ is increasing and continuous in $]0,+\infty[$, with $\lim_{\lambda\to 0^+}\Psi^{-1}(\lambda)=\inf I$.
So, setting $\Psi^{-1}(0)=\inf I$, we think of $\Psi^{-1}$ as an increasing and continuous function in $[0,+\infty[$. Now,
consider the function $\varphi:X\times [0,+\infty[\to {\bf R}$ defined by
$$\varphi(x,\lambda)=\lambda\Phi(x)-J(x)-\int_0^{\lambda}\Psi^{-1}(t)dt$$
for all $(x,\lambda)\in X\times [0,+\infty[$. Of course, for each $x\in X$, the function $\varphi(x,\cdot)$ is concave, while, for each $\lambda>0$,
the function $\varphi(\cdot,\lambda)$ is lower semicontinuous, inf-compact and admits a unique global minimum. Consequently, in view of
Theorem 1.1 of [8], we have
$$\sup_{\lambda>0}\inf_{x\in X}\varphi(x,\lambda)=\inf_{x\in X}\sup_{\lambda>0}\varphi(x,\lambda)\ .$$
But, by continuity, we have
$$\sup_{\lambda\geq 0}\varphi(x,\lambda)=\sup_{\lambda>0}\varphi(x,\lambda)$$
and so
$$\inf_{x\in X}\sup_{\lambda\geq 0}\varphi(x,\lambda)=\inf_{x\in X}\sup_{\lambda>0}\varphi(x,\lambda)\leq
\sup_{\lambda\geq 0}\inf_{x\in X}\varphi(x,\lambda)$$
from which it follows that
$$\sup_{\lambda\geq 0}\inf_{x\in X}\varphi(x,\lambda)=\inf_{x\in X}\sup_{\lambda\geq 0}\varphi(x,\lambda)\ .\eqno{(2.1)}$$
Of course, $\lim_{\lambda\to +\infty}\int_0^{\lambda}\Psi^{-1}(t)dt=+\infty$. By assumption, there is some $x_0\in X$ such that
$\Phi(x_0)=0$, and so we have
$\lim_{\lambda\to +\infty}\varphi(x_0,\lambda)=-\infty$. Hence, $\varphi(x_0,\cdot)$ is sup-compact. Then, from $(2.1)$
it follows that there exists $(\tilde x,\tilde\lambda)\in X\times [0,+\infty[$ such that
$$\sup_{\lambda\geq 0}\varphi(\tilde x,\lambda)=\varphi(\tilde x,\tilde\lambda)=\inf_{x\in X}\varphi(x,\tilde\lambda)\ .$$
Consequently
$$\sup_{\lambda\geq 0}\left(\lambda\Phi(\tilde x)-\int_0^{\lambda}\Psi^{-1}(t)dt\right)=\tilde\lambda\Phi(\tilde x)-\int_0^{\tilde\lambda}\Psi^{-1}(t)dt
\eqno{(2.2)}$$
and
$$\inf_{x\in X}(\tilde\lambda\Phi(x)-J(x))=\tilde\lambda\Phi(\tilde x)-J(\tilde x)\ .\eqno{(2.3)}$$
Notice that $\tilde\lambda>0$. Indeed, otherwise, by $(2.3)$, $\tilde x$ would be a global maximum of $J$, against an assumption.
On the other hand, by $(2.2)$ it follows that $\Psi^{-1}(\tilde\lambda)=\Phi(\tilde x)$.
Hence, $\Phi(\tilde x)\in I$,  $\tilde\lambda=\Psi(\Phi(\tilde x))$, and $(2.3)$ gives
the conclusion, for the existence part. Now, let us prove the uniqueness of $\tilde x$. So, let $\tilde y\in X$ be such that
$\Phi(\tilde y)\in I$ and
$$\Psi(\Phi(\tilde y))\Phi(\tilde y)-J(\tilde y)=\inf_{x\in X}(\Psi(\Phi(\tilde y))\Phi(x)-J(x))\ .$$
Arguing by contradiction, assume that $\tilde y\neq\tilde x$. We claim that $\Phi(\tilde y)=\Phi(\tilde x)$. Indeed, if, for instance,
$\Phi(\tilde y)<\Phi(\tilde x)$, then $\Psi(\Phi(\tilde y))<\Psi(\Phi(\tilde x))$, and so, by Proposition 3.1 of [8], we would have
$\Phi(\tilde y)>\Phi(\tilde x)$. So, if we set $\lambda:=\Phi(\tilde x)$, the points $\tilde x$ and $\tilde y$ would be two distinct global
minima of the function $\lambda\Phi-J$, against an assumption. The proof is complete.

\hfill $\bigtriangleup$\par
\medskip
Here is our main result.\par
\medskip
THEOREM 2.2. - {\it Assume that there exists an open interval $I\subseteq ]0,+\infty[$ such that the restriction of $K$ to $I$ is increasing
and $K(I)=]0,+\infty[$. Let $f:[0,+\infty[\to [0,+\infty[$ be a continuous function, with $f(0)=0$, such that the function
$\xi\to {{f(\xi)}\over {\xi}}$ is decreasing in $]0,+\infty[$ and $\lim_{\xi\to +\infty}{{f(\xi)}\over {\xi}}=0$,
$\lim_{\xi\to 0^+}{{f(\xi)}\over {\xi}}=+\infty\ .$
\par
Then, for each $\alpha\in L^{\infty}(\Omega)$, with $\alpha>0$, the problem
$$\cases {-K\left(\int_{\Omega}|\nabla u(x)|^2dx\right)\Delta u=\alpha(x) f(u)
 & in
$\Omega$\cr & \cr u>0 & in $\Omega$\cr & \cr
u=0 & on
$\partial \Omega$\cr & \cr
\int_{\Omega}|\nabla u(x)|^2dx\in I\cr}$$
has a unique weak solution $\tilde u$ which is the unique global minimum in $H^1_0(\Omega)$ of the functional
$$u\to {{1}\over {2}}K\left(\int_{\Omega}|\nabla \tilde u(x)|^2dx\right)\int_{\Omega}|\nabla u(x)|^2dx-\int_{\Omega}
\left (\alpha(x)\int_0^{u^+(x)}f(t)dt\right)dx\ .$$}\par
\smallskip
PROOF. First of all, extend $f$ to ${\bf R}$ putting $f(\xi)=0$ for all $\xi<0$.
We are going to apply Theorem 2.1 taking $X=H^1_0(\Omega)$ endowed with the weak topology, $\Psi=K$ and defining $\Phi, J$ by
$$\Phi(u)=\int_{\Omega}|\nabla u(x)|^2dx\ ,$$
$$J(u)=2\int_{\Omega}\alpha(x)F(u^+(x))dx$$
for all $u\in H^1_0(\Omega)$, where $F(\xi)=\int_0^{\xi}f(t)dt$. 
The functionals $\Phi, J$
are $C^1$ with derivatives given by
$$\Phi'(u)(v)=2\int_{\Omega}\nabla u(x)\nabla v(x)dx$$
$$J'(u)(v)=2\int_{\Omega}\alpha(x)f(u(x))v(x)dx$$
for all $u,v\in H^1_0(\Omega)$. Moreover, due to the sub-critical growth of $f$, $J$ is
sequentially weakly continuous. Fix $\lambda>0$ and choose $\epsilon\in \left]0,{{\lambda\lambda_1}\over {2}}\right[$,
where $\lambda_1$ is the constant defined in Theorem 1.1.
Since $\lim_{\xi\to +\infty}{{F(\xi)}\over {\xi^{2}}}=0$, there is $c_{\epsilon}>0$ such that
$$F(\xi)\leq {{\epsilon}\over {M}}|\xi|^2+c_{\epsilon}$$
for all $\xi\in {\bf R}$, where $M=\hbox {\rm ess sup}_{\Omega}\alpha$. Consequently, we have
$$J(u)\leq 2\epsilon\int_{\Omega}|u(x)|^2dx+2c_{\epsilon}\int_{\Omega}\alpha(x)dx$$
and so
$$\lambda\Phi(u)-J(u)\geq \lambda\int_{\Omega}|\nabla u(x)|^2dx-2\epsilon\int_{\Omega}|u(x)|^2dx-2c_{\epsilon}\int_{\Omega}
\alpha(x)dx$$
$$\geq \left ( \lambda-{{2\epsilon}\over {\lambda_1}}\right)\int_{\Omega}|\nabla u(x)|^2dx-c_{\epsilon}\int_{\Omega}
\alpha(x)dx$$
for all $u\in H^1_0(\Omega)$. Hence, due to the choice of $\epsilon$, we have
$$\lim_{\|u\|\to +\infty}(\lambda\Phi(u)-J(u))=+\infty\ .$$
This fact, jointly with the reflexivity of $H^1_0(\Omega)$ and the Eberlein-Smulyan theorem, implies that the sequentially weakly
lower semicontinuous functional $\lambda\Phi-J$ is weakly inf-compact. We now show that it has a unique global minimum in $H^1_0(\Omega)$.
Indeed, its critical points are exactly the weak solutions of the problem
$$\cases {-\Delta u={{1}\over {\lambda}}\alpha(x)f(u)
 & in
$\Omega$\cr & \cr 
u=0 & on
$\partial \Omega$\ .\cr}$$
In turn, since the right-hand side of the equation is non-negative, the non-zero weak solutions of the problem are
positive in $\Omega$. Moreover, since, for each $x\in \Omega$, 
the function ${{\alpha(x)f(\xi)}\over {\lambda}\xi}$
is decreasing in $]0,+\infty[$, Theorem 1 of [2] ensures that the problem has at most one positive weak solution, and so it has at most one
non-zero weak solution. As a consequence, we infer that the functional $\lambda\Phi-J$
has a unique global minimum in $H^1_0(\Omega)$, since otherwise, in view of Corollary 1 of [6], it would have at least three critical
points. Now, fix any positive function $u\in C^2(\Omega)\cap C^0(\overline {\Omega})$, with $u=0$ on $\partial\Omega$, such
that 
$$\lambda_1\int_{\Omega}|u(x)|^2dx=\int_{\Omega}|\nabla u(x)|^2dx\ .$$
 Also, fix $\gamma\in ]0,\hbox {\rm ess sup}_{\Omega}\alpha[$.
Of course, the set
$$\Omega_{\gamma}:=\{x\in \Omega : \alpha(x)\geq \gamma\}$$
has a positive measure. Furthermore, fix $M>0$ so that
$$M>{{\lambda\lambda_1\int_{\Omega}|u(x)|^2dx}\over {\gamma\int_{\Omega_{\gamma}}|u(x)|^2dx}}\ .$$
Since $\lim_{\xi\to 0^+}{{F(\xi)}\over {\xi^2}}=+\infty$, there is $\delta>0$ such that
$$F(\xi)\geq M\xi^2$$
for all $\xi\in [0,\delta]$. Now, set $v=\mu u$, where $\mu={{\delta}\over {\sup_{\overline {\Omega}}u}}$. Then, we have
$$J(v)\geq 2M\int_{\Omega}\alpha(x)|v(x)|^2dx\geq 2M\gamma\int_{\Omega_{\gamma}}|v(x)|^2>
\lambda\lambda_1\int_{\Omega}|v(x)|^2dx=\lambda\int_{\Omega}|\nabla v(x)|^2dx=\lambda\Phi(v)\ .$$
This shows that $0$ is not a global minimum for the functional $\lambda\Phi-J$. Consequently, the global minimum of this functional
agrees with its only non-zero critical point. Finally, let us show that $J$ has no global maxima. Arguing by contradiction, suppose
that $\hat u\in H^1_0(\Omega)$ is a global maximum of $J$. Clearly, $J(\hat u)>0$. Consequently, the set
$$A:=\{x\in \Omega : \alpha(x)f(\hat u(x))>0\}$$
has a positive measure. Fix a closed set $C\subset A$ of positive measure. Let $v\in H^1_0(\Omega)$ be such that $v\geq 0$ and
$v(x)=1$ for all $x\in C$. Then, we have
$$\int_{\Omega}\alpha(x)f(\hat u(x))v(x)dx\geq \int_{C}\alpha(x)f(\hat u(x))dx>0$$
and so $J'(\hat u)\neq 0$, which is absurd.
Therefore, each assumption of Theorem 2.1 is satisfied. As a consequence, there exists a unique $\tilde u\in H^1_0(\Omega)$,
with $\int_{\Omega}|\nabla \tilde u(x)|^2dx\in I$, such that
$$K\left(\int_{\Omega}|\nabla \tilde u(x)|^2dx\right)\int_{\Omega}|\nabla\tilde u(x)|^2dx-2\int_{\Omega}\alpha(x)F(\tilde u(x))dx$$
$$=\inf_{u\in H^1_0(\Omega)}\left(K\left(\int_{\Omega}|\nabla \tilde u(x)|^2dx\right)\int_{\Omega}|\nabla u(x)|^2dx-
2\int_{\Omega}\alpha(x)F(u^+(x))dx\right)\ .$$
Clearly, from what seen above, the function $\tilde u$ satisfies the conclusion.\hfill $\bigtriangleup$\par
\medskip
REMARK 2.1. - In Theorem 1.1, the inequality satisfied by the solution follows directly from Theorem 1.2 of [7].\par
\medskip
{\bf Acknowledgements.} The author has been supported by the Gruppo Nazionale per l'Analisi Matematica, la Probabilit\`a e 
le loro Applicazioni (GNAMPA) of the Istituto Nazionale di Alta Matematica (INdAM) and by the Universit\`a degli Studi di Catania, PIACERI 2020-2022, Linea di intervento 2, Progetto ”MAFANE”.
\vfill\eject
\centerline {\bf References}\par
\bigskip
\bigskip
\noindent
[1]\hskip 5pt C. O. ALVES and F. J. S. A. CORR\^EA, {\it On existence of solutions for a class of problem
involving a nonlinear operator}, Comm. Appl. Nonlinear Anal., {\bf 8} (2001), 43-56.\par
\smallskip
\noindent
[2]\hskip 5pt H. BREZIS and L. OSWALD, {\it Remarks on sublinear elliptic equations}, Nonlinear Anal., {\bf 10} (1986), 55-64.\par
\smallskip
\noindent
[3]\hskip 5pt G. M. FIGUEIREDO and  A. SU\'AREZ, {\it Some remarks on the comparison principle in Kirchhoff equations}, Rev. Mat. Iberoam.,
 {\bf 34} (2018), 609-620.\par
\smallskip
\noindent
[4]\hskip 5pt G. M. FIGUEIREDO, C. MORALES-RODRIGO, J. R. SANTOS-J\'UNIOR and A. SU\'AREZ, {\it Study of a nonlinear Kirchhoff equation with non-homogeneous material}, J. Math. Anal. Appl., {\bf 416} (2014), 597-608.\par
\smallskip
\noindent
[5]\hskip 5pt P. PUCCI and V. D. R\u ADULESCU, {\it Progress in nonlinear Kirchhoff problems}, Nonlinear Anal., {\bf 186} (2019), 1-5.\par
\smallskip
\noindent
[6]\hskip 5pt  P. PUCCI and J. SERRIN, {\it A mountain pass theorem},
J. Differential Equations, {\bf 60} (1985), 142-149.\par
\smallskip
\noindent
[7]\hskip 5pt B. RICCERI, {\it A new existence and localization theorem
for the Dirichlet problem}, Dynam. Systems Appl.,  {\bf 22}  (2013),   317-324.\par
\smallskip
\noindent
[8]\hskip 5pt B. RICCERI, {\it On a minimax theorem: an improvement, a new proof and an overview of its applications},
Minimax Theory Appl., {\bf 2} (2017), 99-152.\par
\smallskip
\noindent
[9]\hskip 5pt B. RICCERI, {\it Multiplicity theorems involving functions with non-convex range}, Stud. Univ. Babe\c{s}-Bolyai Math., {\bf 68} (2023), 125-137.\par
\smallskip
\noindent
[10]\hskip 5pt T. SHIBATA, {\it Bifurcation diagrams of one-dimensional Kirchhoff-type equations}
Adv. Nonlinear Anal.,  {\bf 12} (2023), 356-368.\par
\smallskip
\noindent
[11]\hskip 5pt W. WANG and W. TANG, {\it Bifurcation of positive solutions for a nonlocal problem}, Mediterr. J. Math., {\bf 13} (2016), 3955-3964.\par

\bigskip
\bigskip
Department of Mathematics and Informatics\par
University of Catania\par
Viale A. Doria 6\par
95125 Catania, Italy\par
{\it e-mail address:} ricceri@dmi.unict.it

\bye